\documentclass[a4paper]{article}
\usepackage{amssymb}
\usepackage{amsthm}
\usepackage{tikz}
\usepackage{mathtools}
\usepackage[title,toc,page,header]{appendix}
\usepackage{titletoc}

\usepackage[mathscr]{euscript} 
\DeclareFontFamily{OT1}{pzc}{}
\DeclareFontShape{OT1}{pzc}{m}{it}{<-> s * [1.10] pzcmi7t}{}
\DeclareMathAlphabet{\mathpzc}{OT1}{pzc}{m}{it} 
\dottedcontents{section}[2.5em]{\bfseries}{2.5em}{1pc}

\usetikzlibrary{calc}
\usetikzlibrary{matrix,arrows,decorations}
\newtheorem{theorem}{Theorem}[section]
\newtheorem{lemma}[theorem]{Lemma}

\newtheorem{remark}[theorem]{Remark}

\newtheorem*{conjecture*}{Conjecture}

\newtheorem*{remark*}{Remark}
\usepackage{hyperref}

\title{Geometric realization and its variants}

\author{Yi-Sheng Wang}
 
\begin{document}
\maketitle
\begin{abstract}
In this paper, we present a unified approach using model category theory and an associative law to compare some classic variants of the geometric realization functor.        
\end{abstract}  
 
\section{Introduction} 
Given an internal category $\mathcal{C}$ in $\mathfrak{Top}$, the category of weakly Hausdorff $k$-spaces, there are at least three different internal categories in $\mathfrak{Top}$ associated to it:
\[\mathcal{C}^{\mathbb{N}},\hspace*{.5em} \mathcal{C}^{\operatorname{fat}},\hspace*{.3em} \text{ and } \hspace*{.3em}\mathcal{C}^{\operatorname{simp}}.\]
The category $\mathcal{C}^{\mathbb{N}}$ is Segal's unraveled category defined as the subcategory of the product category $\mathcal{C}\times \mathbb{N}$ given by deleting the morphisms $(f,i\leq i)$ with $f\neq \operatorname{id}$, where $\mathbb{N}$ is the linearly ordered set of the natural numbers \cite{Se1}. The geometric realization of the nerve $\vert\operatorname{Ner}_{\cdot}\mathcal{C}^{\mathbb{N}}\vert$ generalizes Milnor's classifying space of a topological group \cite{Mil2}, \cite{Mil3}, \cite{Ha}, \cite{Mo1}, \cite{Mo3}.
 
To define the internal category $\mathcal{C}^{\operatorname{fat}}$ in $\mathfrak{Top}$, we let $\mathcal{OR}$ be the one-object category with morphisms consisting of two elements $0,1$ and composition of morphisms given by the truth table for the operator $\mathsf{or}$; meaning, the composition $a\circ b$ is $1$ when $a$ or $b$ is $1$, and otherwise, it is $0$. Then $\mathcal{C}^{\operatorname{fat}}$ is the subcategory of $\mathcal{C}\times \mathcal{OR}$ given by deleting those morphisms $(f,a)$ with $a=0$ and $f\neq \operatorname{id}$. This construction is functorial, and the geometric realization $\vert\operatorname{Ner}_{\cdot}\mathcal{C}^{\operatorname{fat}}\vert$ is canonically homeomorphic to the fat realization $\vert\vert\operatorname{Ner}_{\cdot}\mathcal{C}\vert\vert$. 

The category $\mathcal{C}^{\operatorname{simp}}$ is the simplex category of the nerve $\operatorname{Ner}_{\cdot}\mathcal{C}$ \cite{Se3}. Its spaces of objects and morphisms are the disjoint unions 
\[\coprod_{[n]}\operatorname{Ner}_{n}\mathcal{C}\hspace*{1em}\text{and}\hspace*{1em}\coprod\limits_{\mathclap{[n]\rightarrow [m]}}\operatorname{Ner}_{m}\mathcal{C},\hspace*{1em}\text{respectively}.\]
The geometric realization $\vert\operatorname{Ner}_{\cdot}\mathcal{C}^{\operatorname{simp}}\vert$ computes the homotopy colimit of the simplicial space $\operatorname{Ner}_{\cdot}\mathcal{C}$ \cite{Du1}, \cite{Hi2}.

Each realization has its own advantage and plays a part in the development of topology. The geometric realization is an important construction in algebraic $K$-theory and delooping theory, and also, it connects the category of simplicial sets $s\operatorname{Sets}$ with the category $\mathfrak{Top}$, making combinatorial methods available in topology. On the other hand, the Segal unraveling construction bridges the gaps between geometry and homotopy theory as, given any topological groupoid $\mathcal{G}$, Segal's construction always gives us the classifying space of numerable $\mathcal{G}$-structures. In fact, by the construction, the space $\vert\operatorname{Ner}\mathcal{G}^{\mathbb{N}}\vert$ admits a numerable universal $\mathcal{G}$-structure \cite[Appendix]{Bo}. On the contrary, the geometric realization of the nerve of a topological groupoid does not always give us the right homotopy type of the classifying space, for example, a topological group that is homeomorphic to the Cantor set; we learn this example from A. Henriques on MathOverflow. However, if replacing geometric realization with fat realization, we get the right homotopy type of the classifying space \cite{HG}. The third construction $\mathcal{C}^{\operatorname{simp}}$ is of importance in model category theory as it computes homotopy colimit \cite{Hi2}, with respect to the projective model structure on the category of simplicial spaces. Also, it is a useful tool in proving theorems (e.g. \cite{Wa3}). Comparison theorems between these constructions allow us to choose models appropriate to different problems and help us understand the geometric meanings of homotopy-theoretic constructions---for instance, the geometric meaning of the algebraic $K$-theory space of a structured category, e.g. an exact category or a Waldhausen category internal in $\mathfrak{Top}$.    

These constructions can be generalized to simplicial spaces or even simplicial objects in a topologically (simplicially) enriched model category $\mathcal{M}$. Furthermore, given a simplicial object $X_{\cdot}$ in $\mathcal{M}$, they are connected by the natural morphisms
\begin{equation}\label{Intro:fivespaces}
X^{\mathbb{N}}_{\cdot}\xrightarrow{\pi} X^{\operatorname{fat}}_{\cdot}\xrightarrow{q} X_{\cdot}\xleftarrow{\mathfrak{l}} X_{\cdot}^{\operatorname{simp}}.
\end{equation} 
The map $\mathfrak{l}$ is the Bousfield-Kan map (or the last vertex map); in the case of simplicial spaces, it has been studied in \cite{BK}, \cite{Se3}, and in the case of simplicial objects in $\mathcal{M}$, \cite{Hi2}, \cite{Du1}. The map $q$, the canonical quotient map from the fat realization to the geometric realization, is rather well understood; \cite{Se3} and \cite{tD} treat the case of simplicial spaces and \cite{Du1} the case of simplicial objects in $\mathcal{M}$. On the contrary, the map $\pi$ is less studied, and \cite{tD} is the only reference dealing with the map $\pi$ that we can find in the literature. \cite[Proposition $2$]{tD} asserts that the map $\pi$ is a homotopy equivalence in the case of simplicial spaces, but it appears that the proof contains some gaps (see Remark \ref{themaprho}); nevertheless, it remains a very promising assertion. In this paper, we present a unified approach, due to Segal, that allows us to compare these constructions simultaneously; the approach is based on an associative law implicitly used in \cite[Appendix $A$]{Se3} and the Reedy model structure on $s\mathcal{M}$, the category of simplicial objects in $\mathcal{M}$. With this approach, we obtain two comparison theorems that recover and generalize most of the known results we know of concerning the relation between these four realization functors, and in particular, we obtain a complete proof of a generalized tom Dieck theorem.    

Our approach relies heavily on a generalized Segal lemma (\cite[Lemma $A.5$]{Se3}) for $\mathfrak{Top}$-enriched model categories.
\begin{theorem}\label{Intro:theSegaltheorem}
Let $s\mathcal{M}$ and $c\mathfrak{Top}$ be the Reedy model categories of simplicial objects in $\mathcal{M}$ and cosimplicial spaces, respectively, and suppose the morphisms $f_{\cdot}:X_{\cdot}\rightarrow Y_{\cdot}\in s\mathcal{M}$ and $g^{\cdot}:I^{\cdot}\rightarrow J^{\cdot}\in c\mathfrak{Top}$ are level-wise weak equivalences between cofibrant objects. Then the induced map between the associated coends
\[\mathclap{\int^{\triangle}}f_{\cdot} \square g^{\cdot}\hspace*{-.1em}:  \hspace*{1em}\mathclap{\int^{\triangle}} X_{\cdot} \square I^{\cdot} \rightarrow\hspace*{.5em} \mathclap{\int^{\triangle}} Y_{\cdot}\square J^{\cdot} \]
is a weak equivalence in $\mathcal{M}$, where $\square$ is the tensor product functor from $\mathcal{M}\times\mathfrak{Top}$ to $\mathcal{M}$, and \hspace*{.7em}$\mathclap{\int^{\triangle}}\hspace*{.5em}(-)\square(-)\in \mathcal{M}$ denotes the coend of a simplicial object in $\mathcal{M}$ and a cosimplicial space. $\triangle$ is the simplex category. 
\end{theorem}    
The category $\mathfrak{Top}$ admits at least three different model structures \cite{MP} and the theorem applies to all of them.
However, for our purpose, we are primarily concerned with the Str\o m model structure. The simplicial version of Theorem \ref{Intro:theSegaltheorem} is discussed in details in \cite[$18.4$]{Hi2}; our proof is different from \cite{Hi2} and based on the decomposition of latching objects in \cite[VII]{GJ}. $\mathfrak{Top}$-($s\operatorname{Sets}$-)enriched model categories of interest to us are the category $\mathpzc{Top}$ \cite{MP}, the category of simplicial sets $s\operatorname{Sets}$ \cite{GJ}, the category of (simplicial) spectra or $\Gamma$-spaces \cite{BF}, and the category of chain complexes \cite{MP}.

\subsection{Main Theorems}  
\begin{theorem}\label{Intro:Thm1}
For any cofibrant object $X_{\cdot}$ in $s\mathcal{M}$, the natural morphisms 
\[\hspace*{.5em}\mathclap{\int^{\triangle}} X^{\mathbb{N}}_{\cdot}\square \triangle^{\cdot}  \xrightarrow{\pi}\hspace*{.5em}\mathclap{\int^{\triangle}} X^{\operatorname{fat}}_{\cdot}\square \triangle^{\cdot}  \xrightarrow{q}\hspace*{.5em} \mathclap{\int^{\triangle}} X_{\cdot} \square \triangle^{\cdot} \xleftarrow{\mathfrak{l}}\hspace*{.5em} \mathclap{\int^{\triangle}} X_{\cdot}^{\operatorname{simp}} \square \triangle^{\cdot}\]
are weak equivalences in $\mathcal{M}$, where the cosimplicial space $\triangle^{\cdot}$ is given by the geometric realization of the standard $n$-simplex $\triangle^{n}_{\cdot}$ in $s\operatorname{Sets}$.  
\end{theorem}
The theorem implies that the geometric realizations of the simplicial spaces in \eqref{Intro:fivespaces} are homotopy equivalent when $X_{\cdot}$ is a proper simplicial space.  
\begin{theorem}\label{Intro:Thm2}
Given a level-wise cofibrant object $X_{\cdot}$ in $s\mathcal{M}$, the following coends 
\[\mathclap{\int^{\triangle}} X_{\cdot}^{\mathbb{N}}\square \triangle^{\cdot} \xrightarrow[\simeq]{\pi} \hspace*{.5em} \mathclap{\int^{\triangle}} X_{\cdot}^{\operatorname{fat}}\square  \triangle^{\cdot} \simeq\hspace*{.5em} \mathclap{\int^{\triangle}} X_{\cdot}^{\operatorname{simp}} \square \triangle^{\cdot}\]
are weakly homotopy equivalent---connected by a zig-zag of weak equivalences.
\end{theorem}
Applying the theorem to the Str\o m model category $\mathpzc{Top}$, we see the projection 
\begin{equation}\label{tomDieckobservation}
\pi:\vert\vert X_{\cdot}\times S_{\cdot}\vert\vert=\vert X_{\cdot}^{\mathbb{N}}\vert\rightarrow \vert X_{\cdot}^{\operatorname{fat}}\vert=\vert\vert X_{\cdot}\vert\vert
\end{equation} 
is a homotopy equivalence, for every simplicial space $X_{\cdot}$, and hence recover \cite[Proposition $2$]{tD}, where $S_{\cdot}$ is the semi-simplicial set defined by $S_{n}:=\{i_{0}<...<i_{n}\mid i_{j}\in\mathbb{N}\}$. The idea of the proof comes from \cite[p309-310]{Se3}, where a kind of associativity is implicitly employed---we interpret it as an associative law in infinite-dimensional linear algebra, namely
\[(v^{T}A)w=v^{T}(Aw),\] 
for any $\infty$-by-$\infty$ matrix $A$ and column vectors $v$ and $w$.

In the last section, we define a map $\tau:  \int^{\triangle} X_{\cdot}^{\operatorname{fat}}\square  \triangle^{\cdot}\rightarrow \int^{\triangle} X_{\cdot}^{\mathbb{N}}\square \triangle^{\cdot}$ to replace the map $\rho$ constructed in \cite[p.47]{tD}\footnote{The construction of $\rho$ appears not to give a well-defined map (see Remark \ref{themaprho}).} and prove that, under the same condition of Theorem \ref{Intro:Thm2}, the map $\tau$ is a homotopy inverse to the map $\pi$.       

The author wishes to thank Sebastian Goette for suggesting the construction of the map $\tau$. He gratefully acknowledges use of facilities at and the financial support from Mathematical Research Institute of Oberwolfach. He thanks God, who gives him life and sustains him.   
  
\section{Left Kan extension}
Let $\triangle_{+}$ denote the subcategory of the simplex category $\triangle$ consisting of injective morphisms. Then a semi-simplicial object in $\mathcal{M}$ is a functor $X_{\cdot}:\triangle_{+}^{\operatorname{op}}\rightarrow \mathcal{M}$ and its left Kan extension $\mathfrak{L}X_{\cdot}$, with respect to the inclusion $\triangle_{+}\hookrightarrow  \triangle$, is given by
\[\mathfrak{L}X_{n}:=\coprod_{\mathclap{[n]\overset{v}{\twoheadrightarrow }[k]}}X_{k},\] 
where $\twoheadrightarrow$ (resp. $\rightarrowtail$) stands for a surjective (resp. injective) morphism \cite[Chapter $X$]{Mac2}, \cite[p.42]{tD}. The simplicial structure of $\mathfrak{L}X_{\cdot}$ can be described as follows: Given a morphism $[n^{\prime}]\xrightarrow{u} [n]$, we let $[n^{\prime}]\overset{s_{u,v}}{\twoheadrightarrow} [k^{\prime}]\overset{i_{u,v}}\rightarrowtail [k]$ be the unique factorization of the composition $[n^{\prime}]\xrightarrow{u} [n]\overset{v}{\twoheadrightarrow} [k]$. Then $u^{\ast}:\mathfrak{L}X_{n}\rightarrow \mathfrak{L}X_{n^{\prime}}$ is given by 
\[\coprod_{\mathclap{[n]\overset{v}{\twoheadrightarrow} [k]}}X_{k}\xrightarrow{\coprod\limits_{v}i_{u,v}^{\ast}} \coprod_{\mathclap{[n^{\prime}]\overset{s_{u,v}}{\twoheadrightarrow} [k^{\prime}]}}X_{k^{\prime}}.\]   
The following generalizes \cite[Lemma $1$]{tD}.
\begin{lemma}\label{GentomDieckLemma}
Given a semi-simplicial object $X_{\cdot}$, there is a canonical isomorphism
\[\mathclap{\int^{\triangle_{+}}}\hspace*{.5em} X_{\cdot}\square \triangle^{\cdot} \xrightarrow{\cong}\hspace*{.5em} \mathclap{\int^{\triangle}}\hspace*{.5em} \mathfrak{L}X_{\cdot}\square \triangle^{\cdot}.\]
\end{lemma} 
\begin{proof}
The isomorphism is given by the inclusion
\[ \coprod_{[n]} X_{n}\square \triangle^{n}\xrightarrow{\coprod\operatorname{id}}  \coprod_{\mathclap{\substack{[n]\\ [n]\xrightarrow{=}[n]}}} X_{n}\square \triangle^{n}\subset  \coprod_{\mathclap{\substack{[n]\\ [n]\twoheadrightarrow [m]}}} X_{m}\square \triangle^{n}.\] 
To define its inverse, we observe that the following morphisms
\begin{align*}
\coprod_{\mathclap{\substack{[n]\\ [n]\overset{v}{\twoheadrightarrow} [s]}}}X_{s}\square \triangle^{n}& \xrightarrow{\coprod\limits_{\tiny \mathclap{ [n]; v}  }\operatorname{id}\square v_{\ast}} \coprod_{[s]} X_{s}\square \triangle^{s};\\
\coprod_{\mathclap{\substack{[n]\xrightarrow{u} [m]\\ [m]\overset{v}{\twoheadrightarrow} [s]}}}X_{s}\square \triangle^{n}&\xrightarrow{ \coprod\limits_{\mathclap{\tiny u,v}}\operatorname{id}\square s_{u,v,\ast}} \coprod_{\mathclap{[s^{\prime}]\overset{i_{u,v}}{\rightarrowtail}[s]}} X_{s}\square \triangle^{s^{\prime}}\\
\end{align*}
respect face and degeneracy maps in $\mathfrak{L}X_{\cdot}$ and face maps in $X_{\cdot}$, and hence they induce a morphism   
\[  \mathclap{\int^{\triangle}}\hspace*{.5em}\mathfrak{L}X_{\cdot}\square\triangle^{\cdot}   \rightarrow \hspace*{.7em} \mathclap{\int^{\triangle_{+}}}\hspace*{.5em}X_{\cdot}\square \triangle^{\cdot}.\]
\end{proof}

\noindent
\textbf{Example $1$:}
Let $X_{\cdot}$ be a simplicial object in $\mathcal{M}$. Regarding it as a semi-simplicial object by the inclusion $\triangle_{+}\hookrightarrow \triangle$, we denote its left Kan extension by $X_{\cdot}^{\operatorname{fat}}$, and there is a canonical projection $X_{\cdot}^{\operatorname{fat}}\rightarrow X_{\cdot}$ given by the assignment
\[\coprod_{\mathclap{[n]\overset{u}{\twoheadrightarrow} [m]}}X_{m}\xrightarrow{\coprod\limits_{\tiny u}u^{\ast}} \coprod_{n}X_{n}.\]

\noindent
\textbf{Example $2$:}
Given an object $X_{\cdot}$ in $s\mathcal{M}$, we denote the left Kan extension of the semi-simplicial object $X_{\cdot}\times S_{\cdot}$ \eqref{tomDieckobservation} by $X_{\cdot}^{\mathbb{N}}$, whose $n$-th component $X_{n}^{\mathbb{N}}$ can be described as follows: 
\[\coprod_{\mathclap{[n]\twoheadrightarrow [k]\rightarrowtail \mathbb{N}}}X_{k}.\]
There is a canonical projection from $X_{\cdot}^{\mathbb{N}}\rightarrow X_{\cdot}^{\operatorname{fat}}$ given by
\[\coprod_{\mathclap{[n]\twoheadrightarrow [k]\rightarrowtail \mathbb{N}}}X_{k}\xrightarrow{\coprod\operatorname{id}}  \coprod_{\mathclap{[n]\twoheadrightarrow [k]}}X_{k}\]

The constructions $(-)^{\mathbb{N}}$ and $(-)^{\operatorname{fat}}$ can be viewed as functors from $s\mathcal{M}$ to itself, and they generalize Segal's unraveling construction and the fat construction defined in the introduction; namely, the following diagrams are commutative  
%%%%%%%%Two diagrams parallel to each other
\begin{center}
\begin{equation}\label{unravellingconstr}
\begin{tikzpicture}[baseline=(current bounding box.center)]
\node(Lu) at (0,2) {$\operatorname{Cat}^{\mathcal{M}}$};
\node(Ll) at (0,0) {$s\mathcal{M}$}; 
\node(Ru) at (2,2) {$\operatorname{Cat}^{\mathcal{M}}$};
\node(Rl) at (2,0) {$s\mathcal{M}$};

\path[->, font=\scriptsize,>=angle 90] 

(Lu) edge node [above]{$(-)^{\mathbb{N}}$}(Ru)  
(Lu) edge node [right]{$\operatorname{Ner}_{\cdot}$}(Ll)
(Ll) edge node [above]{$(-)^{\mathbb{N}}$}(Rl) 
(Ru) edge node [right]{$\operatorname{Ner}_{\cdot}$}(Rl);

%%%%%%%%%%

\node(Luf) at (4,2) {$\operatorname{Cat}^{\mathcal{M}}$};
\node(Llf) at (4,0) {$s\mathcal{M}$}; 
\node(Ruf) at (6,2) {$\operatorname{Cat}^{\mathcal{M}}$};
\node(Rlf) at (6,0) {$s\mathcal{M}$};

\path[->, font=\scriptsize,>=angle 90] 

(Luf) edge node [above]{$(-)^{\operatorname{fat}}$}(Ruf)  
(Luf) edge node [right]{$\operatorname{Ner}_{\cdot}$}(Llf)
(Llf) edge node [above]{$(-)^{\operatorname{fat}}$}(Rlf) 
(Ruf) edge node [right]{$\operatorname{Ner}_{\cdot}$}(Rlf);
\end{tikzpicture}  
\end{equation}
\end{center}
where $\operatorname{Cat}^{\mathcal{M}}$ is the category of internal categories in $\mathcal{M}$.   
%%Functoriality of other two constructions
    
%%Include the other two cases&  
\begin{lemma}\label{lemmaofhtyequivalences}
The canonical maps of simplicial sets \[\triangle^{n,\mathbb{N}}_{\cdot}\xrightarrow{\pi}\triangle^{n,\operatorname{fat}}_{\cdot}\xrightarrow{q}\triangle^{n}_{\cdot}\xleftarrow{\mathfrak{l}} \triangle^{simp}_{\cdot}\] 
induce homotopy equivalences
\[\triangle^{n,\mathbb{N}}\xrightarrow{\pi} \triangle^{n,\operatorname{fat}}\xrightarrow{q} \triangle^{n}\xleftarrow{\mathfrak{l}} \triangle^{simp} \]
after geometric realization.
\end{lemma} 
\begin{proof}
Since the standard $n$-simplex $\triangle^{n}_{\cdot}$ is the nerve of the linearly ordered set $[n]:=\{0\leq 1...\leq n\}$, it suffices to show that the following functors
\[[n]^{\mathbb{N}}\xrightarrow{\pi}[n]^{\operatorname{fat}}\xrightarrow{q}[n]\xleftarrow{\mathfrak{l}}[n]^{\operatorname{simp}}\] 
induce homotopy equivalences, where $[n]^{\operatorname{simp}}$ is the over category $\triangle\downarrow [n]$. It is clear that $[n]^{\operatorname{simp}}$ have the terminal object $n = n$, and the object $n$ is the terminal object in $[n]^{\operatorname{fat}}$ and $[n]$. Hence, if we can show that the composition $q\circ \pi:[n]^{\mathbb{N}}\rightarrow [n]$ induces a homotopy equivalence, then the lemma follows. To see this, we define an intermediate category $[n]^{\mathbb{N},\prime}$ of $[n]^{\mathbb{N}}$ which consists of objects $(k,l)$ with $k\leq l$ and observe that $q\circ\pi$ can be decomposed into two natural projections 
\begin{align*}
\pi_{1}:[n]^{\mathbb{N}}&\rightarrow [n]^{\mathbb{N},\prime},\\
            (k,l)&\mapsto (k,k)& k\geq l,\\
            (k,l)&\mapsto (k,l)& k\leq l;\\                     
\pi_{2}:[n]^{\mathbb{N},\prime}&\rightarrow [n]\\
                            (k,l)&\mapsto k.
\end{align*}
It is clear that the slice category $\pi_{1}\downarrow (k,l)$ has a terminal object $(k,l)$ and the slice category $\pi_{2}\downarrow k$ has an initial object $(0,0)$. By Quillen's theorem $A$, the functors $\pi_{1}$, $\pi_{2}$ induce homotopy equivalences, and hence the proof is complete.
 
\end{proof} 
 
%%The homotopy equivalence
\section{An associative law} 
\begin{lemma}\label{associativitylemma}
There are canonical isomorphisms
\begin{align*}
\mathclap{\int^{\triangle}}X_{\cdot} \square \triangle^{\cdot,\mathbb{N}} &\xrightarrow{\cong}\hspace*{.5em}\mathclap{\int^{\triangle}}X^{\mathbb{N}}_{\cdot} \square \triangle^{\cdot} \\
\mathclap{\int^{\triangle}} X_{\cdot}\square \triangle^{\cdot,\operatorname{fat}} &\xrightarrow{\cong}\hspace*{.5em}\mathclap{\int^{\triangle}}  X_{\cdot}^{\operatorname{fat}}\square\triangle^{\cdot}\\
\mathclap{\int^{\triangle}}X_{\cdot}  \square \triangle^{\cdot,\operatorname{simp}}  &\xrightarrow{\cong}\hspace*{.5em}\mathclap{\int^{\triangle}} X_{\cdot}^{\operatorname{simp}} \square \triangle^{\cdot} 
\end{align*}
\end{lemma}
In the case of simplicial spaces, the last two isomorphisms has been implicitly used in \cite[p.309]{Se3} and \cite[p.359]{Wa3}, and a detailed explanation of the second isomorphism using the universal property of Kan extension is given in \cite[Lemma $1$]{tD}. Here, we present a unified approach to such isomorphisms, viewing them as a consequence of an associative law in infinite dimensional linear algebra.
\begin{proof}
Firstly, we observe that the assignments 
\begin{align*}
X_{n} \square \triangle^{n,\mathbb{N}}_{k} =\coprod_{\mathclap{\substack{[r]\xrightarrow{u} [n]\\ [k]\twoheadrightarrow [r]\hookrightarrow \mathbb{N}}}}X_{n}&\xrightarrow{\coprod u^{\ast}}  \coprod_{\mathclap{[k]\twoheadrightarrow [r]\hookrightarrow \mathbb{N}}}X_{r}=X^{\mathbb{N}}_{k},\\
X_{n} \square \triangle^{n,\operatorname{fat}}_{k} = \coprod_{\mathclap{\substack{[r]\xrightarrow{u} [n] \\ [k]\twoheadrightarrow [r]}}}X_{n}&\xrightarrow{\coprod u^{\ast}}  \coprod_{\mathclap{[k]\twoheadrightarrow [r]}}X_{r}=X^{\operatorname{fat}}_{k},\\
X_{n}  \square  \triangle^{n,\operatorname{simp}}_{k} =\coprod_{\mathclap{\substack{\tiny [r_{0}]\rightarrow...\\ \rightarrow [r_{k}]\xrightarrow{u} [n]}}} X_{n}&\xrightarrow{\coprod u^{\ast}} \coprod_{\mathclap{\tiny [r_{0}]\rightarrow... \rightarrow [r_{k}]}}X_{r_{k}}=: X^{\operatorname{simp}}_{k}
\end{align*}
induce the isomorphisms, whose inverses are given by the obvious inclusions,
\begin{align}
\mathclap{\int^{\triangle}}X_{\cdot} \square \triangle^{\cdot,\mathbb{N}}_{\cdot} &\xrightarrow{\cong} X^{\mathbb{N}}_{\cdot}\nonumber\\
\mathclap{\int^{\triangle}} X_{\cdot}\square \triangle^{\cdot,\operatorname{fat}}_{\cdot}&\xrightarrow{\cong} X_{\cdot}^{\operatorname{fat}}\label{eq:pf:associativity:Nconstr}\\
\mathclap{\int^{\triangle}} X_{\cdot}  \square \triangle^{\cdot,\operatorname{simp}}_{\cdot} &\xrightarrow{\cong}  X_{\cdot}^{\operatorname{simp}}.\nonumber
\end{align}
Then, as different ways of computing colimits yield the same result, there is an isomorphism
\begin{equation}\label{eq:associativity:colimit}
\mathclap{\int^{\triangle}}\hspace*{.4em} X_{\cdot}\square(\hspace*{.7em}\mathclap{\int^{\triangle}}\hspace*{.5em}\triangle^{\cdot,-}_{\cdot}\square\triangle^{\cdot})\cong \hspace*{.5em}\mathclap{\int^{\triangle}}\hspace*{.5em} (\hspace*{.7em}\mathclap{\int^{\triangle}}\hspace*{.5em}X_{\cdot}\square \triangle^{\cdot,-}_{\cdot} )\square\triangle^{\cdot},
\end{equation}
where $-$ can be $\mathbb{N}$, $\operatorname{fat}$ or $\operatorname{simp}$. By \eqref{eq:pf:associativity:Nconstr}, isomorphism \eqref{eq:associativity:colimit} gives the isomorphism  
\[\mathclap{\int^{\triangle}}  X_{\cdot}\square \triangle^{\cdot,-}\cong \hspace*{.5em}\mathclap{\int^{\triangle}} X_{\cdot}^{-} \square \triangle^{\cdot}.\] 
\end{proof}
\begin{remark}
If we view $(\triangle^{n,\mathbb{N}}_{k})$, $(\triangle^{n,\operatorname{fat}}_{k})$, and $(\triangle^{n,\operatorname{simp}}_{k})$ as matrices and $(X_{n})$ and $(\triangle^{k})$ column vectors, then \eqref{eq:associativity:colimit} resembles an associative law in linear algebra.   
\end{remark}
\section{Comparison theorems} 
This section discuss a generalized version of Segal's lemma \cite[Lemma $A.5$]{Se3} for a $\mathpzc{Top}$-enriched model category $\mathcal{M}$; an analogous version for simplicially enriched model
categories can be found in \cite[Corollary $19.4.13$-$14$]{Hi2}.
\begin{lemma}\label{theSegallemma} 
Let $f_{\cdot}:X_{\cdot}\rightarrow Y_{\cdot}$ and $g^{\cdot}:I^{\cdot}\rightarrow J^{\cdot}$ be level-weak equivalences between cofibrant objects in $s\mathcal{M}$ and $c\mathpzc{Top}$, respectively. Then the induced map between the associated coends
\[\mathclap{\int^{\triangle}} f_{\cdot}\square g^{\cdot}\hspace*{-.1em}:\hspace*{1em}\mathclap{\int^{\triangle}} X_{\cdot}\square I^{\cdot} \rightarrow\hspace*{.5em} \mathclap{\int^{\triangle}} Y_{\cdot}\square  J^{\cdot} \] 
is a weak equivalence in $\mathcal{M}$.  
\end{lemma}
\begin{proof}
The idea has been sketched in \cite[Appendix]{Se3} (see also \cite[p.43]{tD}, \cite[p.11]{Du1}, \cite[p.375]{GJ}). For the sake of completeness, we give a detailed proof here. 
Firstly, we claim that the latching object of any cofibrant object in $s\mathcal{M}$ (resp. $c\mathpzc{Top}$) is cofibrant. Following \cite[p.362-6]{GJ}, we consider the category $\mathcal{O}_{n}$ whose objects are surjective morphisms $[n]\twoheadrightarrow [m]$ with $n>m$ and morphisms from $[n]\twoheadrightarrow [m]$ to $[n]\twoheadrightarrow [m^{\prime}]$ are those morphisms $[m]\rightarrow [m^{\prime}]$ satisfying the commutative diagram
\begin{center}
\begin{tikzpicture}
\node(Ll) at (0,0){$[m]$};
\node(Rl) at (2,0){$[m^{\prime}]$};  
\node(Mu) at (1,1) {$[n]$};

\path[->, font=\scriptsize,>=angle 90] 

(Mu) edge (Rl)  
(Mu) edge (Ll) 
(Ll) edge (Rl);
\end{tikzpicture}
\end{center}
The $n$-th latching object of a (co)simplicial object $X_{\cdot}$, denoted by $L_{n}X_{\cdot}$, is then given by the colimit
\[\operatorname*{colim}\limits_{\tiny\mathclap{[n]\twoheadrightarrow [m] \in \mathcal{O}^{\operatorname{op}}_{n}}}X_{m}.\] 
Now, for each $1\leq k\leq n$, we can define the subcategories of $\mathcal{O}_{n}$
\begin{align*}
\mathcal{M}_{n,k}&:=\{\phi:[n]\twoheadrightarrow [m]\mid \phi(k)\leq k\}\\
\mathcal{M}(k-1)&:=\{\phi:[n]\twoheadrightarrow [m]\mid \phi(k-1)=\phi(k)\}.
\end{align*}
Let $L_{n,k}X_{\cdot}$ be the colimit 
\[\operatorname*{colimt}\limits_{\tiny\mathclap{[n]\twoheadrightarrow [m]\in\mathcal{M}_{n,k}}}X_{m}.\] 
Then there is a filtration of $L_{n}X_{\cdot}$ given by
\begin{equation}\label{pf:segallemma:filtration}
X_{n-1}=L_{n,1}X_{\cdot}\subset L_{n,2}X_{\cdot}\subset...\subset L_{n,n}X_{\cdot}=L_{n}X_{\cdot}
\end{equation} 
and a pushout diagram
\begin{center}
\begin{equation}\label{pf:segallemma:pushoutforlatching}
\begin{tikzpicture}[baseline=(current bounding box.center)]
\node(Lu) at (0,1.5) {$L_{n-1,k}X_{\cdot}$};
\node(Ll) at (0,0) {$L_{n,k}X_{\cdot}$}; 
\node(Ru) at (5,1.4) {$\operatorname*{colim}\limits_{\tiny\mathclap{[n]\twoheadrightarrow [m]\in\mathcal{M}(k)}}X_{m}=X_{n-1}$};
\node(Rl) at (4,0) {$L_{n,k+1}X_{\cdot}$};

\draw[->] (Lu) to node [above]{\scriptsize $s_{k}^{\ast}$}(3.3,1.5);  
\draw[->] (Lu) to node [right]{\scriptsize $s_{k}^{\ast}$}(Ll);
\draw[->] (Ll) to (Rl); 
\draw[->] (4,1) to (Rl);
\end{tikzpicture}
\end{equation}
\end{center}
where $s_{k}:[n]\twoheadrightarrow [n-1]$ is the degeneracy map with $s_{k}(k)=s_{k}(k+1)$. By induction and pushout diagram \eqref{pf:segallemma:pushoutforlatching}, we see filtration \eqref{pf:segallemma:pushoutforlatching} is a sequence of cofibrations and the object $L_{n,k}X_{\cdot}$ is cofibrant, for every $n,k$. In particular, the objects $L^{n}I^{\cdot}$, $L^{n}J^{\cdot}$, $L_{n}X_{\cdot}$, and $L_{n}Y_{\cdot}$ are cofibrant, for every $n$.  

Now, since $\mathcal{M}$ is a $\mathpzc{Top}$-enriched category, given a cofibrant object $Z$ in $\mathcal{M}$ and a cofibrant object $A$ in $\mathpzc{Top}$, the functors below
\begin{align*}
A\square -:\mathcal{M}&\mapsto \mathcal{M}\\
-\square Z:\mathpzc{Top}&\mapsto \mathcal{M}
\end{align*}
preserve cofibrations and weak equivalences between cofibrant objects \cite[Lemmas $14.2.9$; $16.4.5$]{MP}. Thus, we have the following weak equivalences between two cospan of cofibrations
\begin{equation}\label{pf:segallemma:pushout1}
(X_{n}\square L^{n}I^{\cdot}\leftarrow L_{n}X_{\cdot}\square L^{n}I^{\cdot}\rightarrow L_{n}X_{\cdot}\square I^{n})\rightarrow (Y_{n}\square L^{n}J^{\cdot}\leftarrow L_{n}Y_{\cdot}\square L^{n}J^{\cdot}\rightarrow L_{n}Y_{\cdot}\square J^{n}).
\end{equation} 
Since the subcategory of cofibrant objects in a model category is always proper, \eqref{pf:segallemma:pushout1} induces a weak equivalence between the pushouts  
\begin{equation}\label{pf:segallemma:weakeq}
X_{n}\square L^{n}I^{\cdot}\cup_{L_{n}X_{\cdot}\square L^{n}I^{\cdot}}L_{n}X_{\cdot}\square I^{n}\rightarrow Y_{n}\square L^{n}J^{\cdot}\cup_{L_{n}Y_{\cdot}\square L^{n}J^{\cdot}}L_{n}Y_{\cdot}\square J^{n}.
\end{equation} 
Now, recall that the $n$-skeleton object $(\operatorname{sk}_{n}Z_{\cdot})_{\cdot}$ of a simplicial object is defined by first truncating $Z_{\cdot}$ at the $n$-th degree, denoted by $\bar{Z}_{\cdot}$ and then freely throwing the degeneracies \cite[p.354]{GJ}, namely
\[[m]\mapsto (\operatorname{sk}_{n}Z_{\cdot})_{m}:=\operatorname*{colim}\limits_{\tiny \mathclap{\substack{[m]\twoheadrightarrow [k]\\ k\leq n}}}\bar{Z}_{k}.\]
By the definition, we have $\operatorname*{colim}\limits_{n}(\operatorname{sk}_{n}Z_{\cdot})_{\cdot}=Z_{\cdot}$.
Since the $n$-skeleton \hspace*{.7em}$\mathclap{\int^{\triangle}}\hspace*{.5em}(\operatorname{sk}_{n}X_{\cdot})_{\cdot}\square I^{\cdot}$ (resp. \hspace*{.7em}$\mathclap{\int^{\triangle}}\hspace*{.5em}(\operatorname{sk}_{n}Y_{\cdot})_{\cdot}\square J^{\cdot}$) is the pushout of the cospan
\[ \mathclap{\int^{\triangle}}\hspace*{.5em}(\operatorname{sk}_{n-1}X_{\cdot})_{\cdot}\square I^{\cdot}\leftarrow X_{n}\square L^{n}I^{\cdot}\cup_{L_{n}X_{\cdot}\square L^{n}I^{\cdot}}L_{n}X_{\cdot}\square I^{n}\rightarrow X_{n}\square I^{n}\] 
\[(\text{resp. }\hspace*{.7em}\mathclap{\int^{\triangle}}\hspace*{.5em}(\operatorname{sk}_{n-1}Y_{\cdot})_{\cdot}\square J^{\cdot}\leftarrow Y_{n}\square L^{n}J^{\cdot}\cup_{L_{n}Y_{\cdot}\square L^{n}J^{\cdot}}L_{n}Y_{\cdot}\square J^{n}\rightarrow Y_{n}\square J^{n}),\]
and the second arrow in each span is a cofibration---$\mathcal{M}$ is $\mathfrak{Top}$-enriched \cite[Lemma $16.4.5$]{MP}, by induction, we get the weak equivalence   
\[\mathclap{\int^{\triangle}}\hspace*{.5em}\operatorname{sk}_{n}X_{\cdot}\square I^{\cdot}\rightarrow \hspace*{.5em}\mathclap{\int^{\triangle}}\hspace*{.5em}\operatorname{sk}_{n}Y_{\cdot}\square J^{\cdot},\]
for every $n$. The theorem then follows from the fact that the functor \hspace*{.7em}$\mathclap{\int^{\triangle}}\hspace*{.5em}(-)\square I^{\cdot}$ (resp. \hspace*{.7em}$\mathclap{\int^{\triangle}}\hspace*{.5em}(-)\square J^{\cdot}$) is a left adjoint and commutes with colimits.    
\end{proof}
As a corollary of Lemmas \ref{lemmaofhtyequivalences}, \ref{associativitylemma}, and \ref{theSegallemma}, we have the following theorem (compare with \cite[Theorem 18.7.4]{Hi2}, \cite[Proposition $1$]{Se3}, \cite[Proposition $1$]{tD}).  
\begin{theorem}\label{mainthm1}
Given a cofibrant object $X_{\cdot}$ in $s\mathcal{M}$, the canonical maps
\[ \mathclap{\int^{\triangle}}\hspace*{.5em}X_{\cdot}^{\mathbb{N}} \square\triangle^{\cdot}\xrightarrow{\pi}\hspace*{.5em}\mathclap{\int^{\triangle}}\hspace*{.5em}X_{\cdot}^{\operatorname{fat}}\square\triangle^{\cdot}\xrightarrow{q}\hspace*{.5em} \mathclap{\int^{\triangle}}\hspace*{.5em}X_{\cdot} \square\triangle^{\cdot}\xleftarrow{\mathfrak{l}} \hspace*{.5em} \mathclap{\int^{\triangle}}\hspace*{.5em} X_{\cdot}^{\operatorname{simp}}\square\triangle^{\cdot}\]
are weak homotopy equivalences.
\end{theorem}   
The following can be deduced from Theorem \ref{mainthm1}.
\begin{theorem}\label{mainthm2}
Let $X_{\cdot}$ be a level-wise cofibrant object in $s\mathcal{M}$. Then the objects below are weakly equivalent 
\[\mathclap{\int^{\triangle}}\hspace*{.5em}X_{\cdot}^{\mathbb{N}}\square\triangle^{\cdot}\xrightarrow[\simeq]{\pi}\hspace*{.5em}\mathclap{\int^{\triangle}}\hspace*{.5em}X_{\cdot}^{\operatorname{fat}}\square\triangle^{\cdot}\simeq  \hspace*{.5em} \mathclap{\int^{\triangle}}\hspace*{.5em} X_{\cdot}^{\operatorname{simp}} \square\triangle^{\cdot}.\]
\end{theorem} 
\begin{proof}
Let $Y_{\cdot}$ be a cofibrant replacement of $X_{\cdot}$. Then the theorem ensues from the following commutative diagram of weak equivalences
\begin{center}
\begin{tikzpicture}
\node(Lu) at (0,1.5) {$ \int^{\triangle}Y_{\cdot}^{\mathbb{N}} \square\triangle^{\cdot}$};
\node(Ll) at (0,0) {$ \int^{\triangle}X_{\cdot}^{\mathbb{N}}\square\triangle^{\cdot}$};
\node(Mu1) at (2.7,1.5) {$ \int^{\triangle}Y_{\cdot}^{\operatorname{fat}}\square\triangle^{\cdot}$};
\node(Ml1) at (2.7,0) {$ \int^{\triangle}X_{\cdot}^{\operatorname{fat}}\square\triangle^{\cdot}$};
\node(Mu2) at (5.5,1.5) {$ \int^{\triangle}Y_{\cdot}\square\triangle^{\cdot}$}; 
\node(Ru) at (9,1.5) {$ \int^{\triangle} Y_{\cdot}^{\operatorname{simp}} \square\triangle^{\cdot}$};
\node(Rl) at (9,0) {$\int^{\triangle} X_{\cdot}^{\operatorname{simp}} \square\triangle^{\cdot}$};

\path[->, font=\scriptsize,>=angle 90]

(Ll) edge (Ml1)
 
(Lu) edge (Mu1)
(Mu1) edge (Mu2)  
(Lu) edge (Ll)
(Mu1) edge (Ml1)
   
(Ru) edge (Mu2)
 
(Ru) edge (Rl);
\end{tikzpicture}
\end{center}
The horizontal arrows above are weak equivalences by Theorem \ref{mainthm1}. The vertical arrows are weak equivalences by Lemma \ref{theSegallemma} because $Y^{\mathbb{N}}_{\cdot}\rightarrow X^{\mathbb{N}}_{\cdot}$ (resp. $Y^{\operatorname{fat}}_{\cdot}\rightarrow X^{\operatorname{fat}}_{\cdot}$ and $Y_{\cdot}^{\operatorname{simp}}\rightarrow  X_{\cdot}^{\operatorname{simp}}$) is a level-wise weak equivalence between cofibrant objects in $s\mathcal{M}$---the assumption that $X_{\cdot}$ is level-wise cofibrant is used here.
\end{proof}  
\begin{remark}
The coend \hspace*{.7em}$\mathclap{\int^{\triangle}}\hspace*{.5em} X_{\cdot}^{\operatorname{simp}} \square\triangle^{\cdot}$ computes the homotopy colimit \cite[Part $I.4$]{Du1}, \cite[Chapter $18.1$]{Hi2} of the diagram $X_{\cdot}$, and hence, all these variants of geometric realization in Theorem \ref{mainthm2} compute the homotopy colimit of a level-wise cofibrant object in $s\mathcal{M}$ (compare with \cite[Section $17.4$]{Du1}). However, when viewing the homotopy colimit functor as the left Kan extension of the colimit functor, we do not use the Reedy model structure on $s\mathcal{M}$ but the projective one \cite[Section $11.6$]{Hi2}, \cite[Section $5.8$]{Du1}.
\end{remark}
         
\section{The homotopy inverse to $\pi$}
Let $\operatorname{Sd}_{\cdot}\triangle^{n}$ be the semi-simplicial set given by  
\[\operatorname{Sd}_{k}\triangle^{n}:=\{[l_{0}]\rightarrowtail...\rightarrowtail [l_{k}]\rightarrowtail [n]\mid l_{i}<l_{i+1}\leq n \text{ for }i=0...k-1\}.\]
Then the fat realization of $\operatorname{Sd}_{\cdot}\triangle^{n}$\hspace*{.7em} 
\[\mathclap{\int^{\triangle_{+}}}\operatorname{Sd}_{\cdot}\triangle^{n}\times \triangle^{\cdot}\] 
is the barycentric subdivision of $\triangle^{n}$ and canonically homeomorphic to $\triangle^{n}$.

Now, consider the assignment below  
\begin{align*}
\bar{\tau}^{n}:\coprod_{m,k}\triangle^{n}_{m}\times \operatorname{Sd}_{k}\triangle^{m}\times\triangle^{k}&\rightarrow  \coprod_{k}\triangle^{n}_{k}\times S_{k}\times \triangle^{k}\\
(x,[l_{0}]\rightarrowtail [l_{1}]\rightarrowtail...\rightarrowtail [l_{k}]\rightarrowtail [m],t)&\mapsto (u^{\ast}x,l_{0} < l_{1}  <...<l_{k} ,t), 
\end{align*}
where $u:[k]\rightarrow [m]$ (the last vertex map) is defined by letting $u(i)$ be the image of $l_{i}$ under the composition
\[[l_{i}]\rightarrowtail...\rightarrowtail [l_{k}]\rightarrowtail [m],\]
and observe that the assignment descends to a map of coends
\[\tau^{n}:\triangle^{n,\operatorname{fat}}=\hspace*{.7em}\mathclap{\int^{\triangle_{+}}}\hspace*{.5em}\triangle^{n}_{\cdot}\times\hspace*{.8em}\mathclap{\int^{\triangle_{+}}}\hspace*{.5em}\operatorname{Sd}_{\cdot}\triangle^{\cdot}\times\triangle^{\cdot}\rightarrow\hspace*{.7em}\mathclap{\int^{\triangle_{+}}}(\triangle^{n}_{\cdot}\times S_{\cdot})\times\triangle^{\cdot}=\triangle^{n,\mathbb{N}}, \]
where the second identity follows from Lemma \ref{GentomDieckLemma}. Observe also that the composition 
\[\triangle^{n,\operatorname{fat}}\xrightarrow{\tau^{n}}\triangle^{n,\mathbb{N}}\xrightarrow{\pi}\triangle^{n,\operatorname{fat}}\]
is induced by the standard map from the barycentric subdivision of a simplex to itself, and hence the linear homotopy gives a homotopy of cosimplicial spaces 
\[H:\triangle^{\cdot,\operatorname{fat}}\times I\rightarrow \triangle^{\cdot,\operatorname{fat}}\] 
between the identity and $\pi\circ\tau^{n}$. Now, let $X_{\cdot}$ be a simplicial object in $\mathcal{M}$. Then, by Lemma \ref{associativitylemma}, the cosimplicial map $\tau^{\cdot}$ induces a morphism
\[\tau:\hspace*{.7em}\mathclap{\int^{\triangle}}\hspace{.5em}X_{\cdot}^{\operatorname{fat}}\square\triangle^{\cdot}=\hspace*{.7em}\mathclap{\int^{\triangle}}\hspace{.5em}X_{\cdot}\square\triangle^{\cdot,\operatorname{fat}}\rightarrow\hspace*{.7em}\mathclap{\int^{\triangle}}\hspace{.5em}X_{\cdot}\square\triangle^{\cdot,\mathbb{N}}=\hspace*{.7em}\mathclap{\int^{\triangle}}\hspace{.5em}X_{\cdot}^{\mathbb{N}}\square\triangle^{\cdot},\]  
and the homotopy $H$ shows that the morphism $\tau$ is the homotopy inverse to $\pi$ when $X_{\cdot}$ is level-wise cofibrant (Theorem \ref{mainthm2}). Furthermore, the morphism $\tau$ is functorial with respect to $X_{\cdot}$ as we have the commutative diagram below, for any morphism $X_{\cdot}\rightarrow Y_{\cdot}$ in $s\mathcal{M}$, 
\begin{center}
\begin{tikzpicture}
\node(Lu) at (0,1.5) {$ \int^{\triangle}X_{\cdot}\square\triangle^{\cdot,\operatorname{fat}}$};
\node(Ll) at (0,0) {$ \int^{\triangle} Y_{\cdot}\square\triangle^{\cdot,\operatorname{fat}}$}; 
\node(Ru) at (5,1.5) {$ \int^{\triangle}X_{\cdot}\square\triangle^{\cdot,\mathbb{N}}$};
\node(Rl) at (5,0) {$ \int^{\triangle}Y_{\cdot}\square\triangle^{\cdot,\mathbb{N}}$};

\path[->, font=\scriptsize,>=angle 90] 

(Lu) edge node [above] {$\tau$}(Ru)  
(Lu) edge (Ll)
(Ll) edge node [above] {$\tau$}(Rl) 
(Ru) edge (Rl);
\end{tikzpicture}
\end{center} 
Hence, we have proved the following theorem.
\begin{theorem}
Given a simplicial object $X_{\cdot}$ in $\mathcal{M}$, there is a well-defined morphism
\[\tau:\hspace*{.7em}\mathclap{\int^{\triangle}}\hspace{.5em}X_{\cdot}^{\operatorname{fat}}\square\triangle^{\cdot}\rightarrow\hspace*{.7em}\mathclap{\int^{\triangle}}\hspace{.5em}X_{\cdot}^{\mathbb{N}}\square\triangle^{\cdot}\]
such that the composition $\pi\circ\tau$ is left homotopy to $\operatorname{id}$ and $\tau$ is functorial with respect to $X_{\cdot}$. If, in addition, $X_{\cdot}$ is level-wise cofibrant, then the morphism $\tau$ is a homotopy inverse to the morphism $\pi$.   
\end{theorem} 
\begin{remark}\label{themaprho}
To define a homotopy inverse to the map $\pi$ in the case of simplicial spaces, \cite[p.45-7]{tD} considers the following assignment 
\begin{align}\label{tomDieckmap}
X_{n}\times \triangle^{n}&\rightarrow X_{n}\times S_{n}\times \triangle^{n}\\
(y;t_{0},...,t_{n})&\mapsto (y,1<...<n;s_{1,n}(t_{0},...,t_{n}),...,s_{n,n}(t_{0},...,t_{n}))\nonumber, 
\end{align}  
where 
\[s_{j,n}(t_{0},...,t_{n}):=(j+1)\sum_{E}\max(0,\min_{j\in E}t_{j}-\max_{j\not\in E}t_{j})\]
and $E$ runs through all subsets of $[n]$ with $j+1$ elements---the map $s_{j,n}$ is a kind of folding map sending all the $n$-simplices in the barycentric subdivision of an $n$-simplex to the $n$-simplex. However, the assignment does not respect face maps---the first and second components in assignment \eqref{tomDieckmap} should depend on $(t_{0},...,t_{n})$, the coordinate of points in $\triangle^{n}$. The map $\rho$ is used in cohomology theories of spaces with two topologies \cite[Theorem $7.1$]{Mo1}. 
\end{remark}  
%%tom Dieck's theorem (note he is working in the category compactly generated space)
 
\addcontentsline{toc}{section}{\hspace*{2.2em} References} 

\bibliographystyle{alpha}

\bibliography{Reference}  

\end{document}